\documentclass{article}

\usepackage{amsmath}
\usepackage{amsfonts}
\usepackage{amssymb}
\usepackage{amsthm}
\usepackage[all]{xy}
\input diagxy

\theoremstyle{plain}

\newtheorem{theorem}{Theorem}[section]

\newtheorem{proposition}[theorem]{Proposition}
\newtheorem{corollary}[theorem]{Corollary}

\theoremstyle{definition}

\newtheorem{definition}[theorem]{Definition}

\newtheorem{remark}[theorem]{Remark}

\def\To{\mathrel{\xy \ar@{=>} (200,0) \endxy}}
\def\comp{\raisebox{0.2mm}{\ensuremath{\scriptstyle\circ}}}

\newcommand{\del}{\partial}

\renewcommand{\ker}{\ensuremath{\mathsf{Ker\,}}}

\DeclareMathOperator{\ext}{Ext}

\newcommand{\Bc}{\ensuremath{\mathcal{B}}}
\newcommand{\Ac}{\ensuremath{\mathcal{A}}}
\newcommand{\Ab}{\mathsf{Ab}}
\newcommand{\ab}{\mathrm{Ab}}

\newcommand{\Gal}{\ensuremath{\mathsf{Gal}}}
\newcommand{\Sc}{\ensuremath{\mathcal{S}}}

\newcommand{\N}{\ensuremath{\mathbb{N}}}
\newcommand{\Gb}{\ensuremath{\mathbb{G}}}
\newcommand{\Set}{\ensuremath{\mathsf{Set}}}

\newcommand{\Ext}{\ensuremath{\mathsf{Ext}}}
\newcommand{\Eq}{\ensuremath{\mathsf{Eq}}}
\newcommand{\CentrExt}{\ensuremath{\mathsf{Centr}}}
\newcommand{\centr}{\ensuremath{\mathrm{Centr}}}

\newcommand{\Fun}{\ensuremath{\mathsf{Fun}}}

\newcommand{\PCh}{\ensuremath{\mathsf{PCh}}}

\newcommand{\simp}[1]{\Sc #1}
\newcommand{\simpA}{\simp{\Ac}}

\renewcommand{\hom}{\mathrm{Hom}}
\DeclareMathOperator{\op}{op}
\newcommand{\noproof}{\hfill \qed}

\newbox\pullbackbox
\setbox\pullbackbox=\hbox{\xy 0;<1mm,0mm>: \POS(4,0)\ar@{-} (0,0) \ar@{-} (4,4)
\endxy}
\def\pullback{\copy\pullbackbox}

\newbox\pushoutbox
\setbox\pushoutbox=\hbox{\xy 0;<1mm,0mm>: \POS(0,4)\ar@{-} (0,0) \ar@{-} (4,4)
\endxy}
\def\pushout{\copy\pushoutbox}

\begin{document}

\newdir{>>}{{}*!/9pt/@{|}*!/4.5pt/:(1,-.2)@^{>}*!/4.5pt/:(1,+.2)@_{>}}
\newdir{ >>}{{}*!/9pt/@{|}*!/4.5pt/:(1,-.2)@^{>}*!/4.5pt/:(1,+.2)@_{>}}
\newdir{ >}{{}*!/-5.5pt/@{|}*!/-10pt/:(1,-.2)@^{>}*!/-10pt/:(1,+.2)@_{>}}
\newdir2{ >}{{}*-<1.8pt,0pt>@{ >}}
\newdir{^ (}{{}*!/-10pt/@{>}}
\newdir{_ (}{{}*!/-10pt/@{>}}
\newdir{>}{{}*:(1,-.2)@^{>}*:(1,+.2)@_{>}}
\newdir{<}{{}*:(1,+.2)@^{<}*:(1,-.2)@_{<}}

\title{On the second cohomology group\\
in semi-abelian categories}

\author{Marino Gran\footnote{Marino Gran, Lab.\ Math.\ Pures et Appliqu\'ees, Universit\'e du Littoral C\^ote d'Opale, 50,~Rue F.~Buisson, 62228~Calais, France. Email:~\texttt{gran@lmpa.univ-littoral.fr}}\, and Tim Van der Linden\footnote{Tim Van~der~Linden, Vakgroep Wiskunde, Vrije Universiteit Brussel, Pleinlaan~2, 1050~Brussel, Belgium. Email:~\texttt{tvdlinde@vub.ac.be}}}

\maketitle

\begin{abstract}
\noindent We develop some new aspects of cohomology in the context of semi-abelian categories: we establish a Hochschild-Serre 5-term exact sequence extending the classical one for groups and Lie algebras; we prove that an object is perfect if and only if it admits a universal central extension; we show how the second Barr-Beck cohomology group classifies isomorphism classes of central extensions; we prove a universal coefficient theorem to explain the relationship with homology.

\end{abstract}

\section*{Introduction}\label{Section-Introduction}
The notion of \emph{semi-abelian category} introduced by Janelidze, M\'arki and Tholen in \cite{Janelidze-Marki-Tholen} provides a natural context for a unified treatment of many important homological properties of the categories of groups, rings, Lie algebras, crossed modules, $C^\star$-algebras and compact Hausdorff groups. Several results, which are classical in the category of groups, can be extended to any semi-abelian category: among them, let us mention the {$3 \times 3$ Lemma}, the Snake Lemma \cite{Bourn2001}, the long exact homology sequence associated with a proper chain complex, or the Stallings-Stammbach five term exact sequence associated with a short exact sequence \cite{EverVdL2}.
Semi-abelian categories thus provide an elegant answer to the old problem, first considered by Mac\,Lane \cite{MacLane:Duality}, of finding a suitable list of axioms that would reflect the homological properties of groups and rings in the same way as the axioms of abelian category reflect some particular properties of the categories of abelian groups and modules over a ring. The fundamental advances in categorical algebra that made it possible to solve this problem, and on which the notion of semi-abelian category is built, are the discoveries of the subtle properties of Barr-exactness \cite{Barr} and of Bourn-protomodularity \cite{Bourn1991}. We refer to the introduction of \cite{Janelidze-Marki-Tholen} for a detailed description of the historical developments that led to the theory of semi-abelian categories.

In the present article, we prove some new results in the study of the homology and cohomology in a semi-abelian algebraic category. More precisely, we can establish: 
\begin{itemize}
\item a natural Hochschild-Serre $5$-term exact sequence \cite{Hochschild-Serre}
in cohomology with trivial coefficients (Theorem~\ref{Theorem-Cohomology-Sequence}); 
\item an interpretation of the second cohomology group $H^2(Y,A)$ as the group $\CentrExt (Y,A)$ of isomorphism classes of central extensions of $Y$ by an abelian object $A$, equipped with (a generalization of) the Baer sum (Theorem~\ref{Theorem-H^2-Central-Extensions});
\item a universal coefficient theorem relating homology and cohomology (Theorem~\ref{Theorem-Universal-Coefficients}).
\end{itemize}
Thus we simplify some recent investigations in this direction in the context of crossed modules \cite{Carrasco-Homology} or precrossed modules \cite{AL}, and unify them with the classical theory that exists for groups and Lie algebras. Our approach is based on the work of Fr\"ohlich's school---in particular, of Lue \cite{Lue}---and also on the work of Janelidze and Kelly \cite{Janelidze-Kelly}, who recently discovered some unexpected connections between homological algebra and universal algebra (see also \cite{EG, Janelidze-Kelly:Maltsev}). And indeed, the present work benefits from some new categorical tools designed for understanding the universal algebraic property of centrality \cite{Carboni-Pedicchio-Pirovano, Pedicchio, Janelidze-Pedicchio, BG}. For a different approach to such problems, the reader is also referred to \cite{Bourn-Janelidze:Torsors}.

In the first section we collect the main properties of semi-abelian categories that will be needed throughout the paper. In Section~\ref{Section-Centrality} we prove that two notions of central extension are equivalent: the first one is the notion of extension whose kernel pair is central in the sense of Smith \cite{Smith}, while the second one is the notion of extension whose kernel is a central arrow in the sense of Huq \cite{Huq}. In the third section a useful technical property involving exact sequences and central extensions is proven. In Section~\ref{Section-Perfect-Case}, we prove that an object is perfect if and only if it admits a universal central extension: this extends a classical result due to Fr\"ohlich \cite{Froehlich}. In Section~\ref{Section-Cohomology} we first recall the Stallings-Stammbach sequence and Hopf's formula for the second homology object in a semi-abelian category \cite{EverVdL2}. We then obtain a cohomological version of Hopf's formula as well as the Hochschild-Serre $5$-term exact sequence for cohomology: an extension ${f \colon X \to Y}$ with kernel $K$
induces the exact sequence
\[
 0 \to<250> H^{1} Y \to/{{ >}->}/<250>^{H^{1}f} H^{1} X \to<250> \hom \bigl(\textstyle{\frac{K}{[K,X]}},A\bigr) \to<250> H^{2} Y \to<250>^{H^{2}f} H^{2} X.
\]
In Section~\ref{Section-Second-Cohomology} we prove that the second cohomology group $H^{2} (Y,A)$ is isomorphic to the group $\CentrExt (Y,A)$ of isomorphism classes of central extensions of $Y$ by $A$. In the last section we establish a universal coefficient theorem for cohomology in semi-abelian categories, and we give some applications.

\medskip\noindent\textbf{Acknowledgements.} We are grateful to Dominique Bourn for a precious remark concerning Proposition~\ref{Proposition-pushout}. Many thanks also to the referee for several interesting remarks and suggestions.

\section{Protomodular and semi-abelian categories}\label{Section-Revision}

In this section we recall some basic definitions and properties of protomodular and semi-abelian categories, needed throughout the article. We shall always assume that the category $\Ac$ in which we are working is finitely complete.

\begin{definition}\cite{Bourn1991}
A finitely complete category $\Ac$ is {\em protomodular} if it satisfies the following property: given any commutative diagram 
\[
\xymatrix{{\cdot} \ar[r] \ar[d] & {\cdot} \ar@{.>}[d]\ar[r] & {\cdot} \ar[d]  \\ 
{\cdot} \ar[r] & {\cdot} \ar[r] & {\cdot}}
\]
where the dotted vertical arrow is a split epimorphism, the left-hand square and the whole rectangle are pullbacks, the right-hand square is also a pullback.
\end{definition}

Recall that a finitely complete category $\Ac$ is \emph{regular} if (1) every kernel pair has a coequalizer and (2) regular epimorphisms are stable under pulling back. A~regular category is protomodular if and only if given any commutative diagram as above, where the dotted vertical arrow is a regular epimorphism, the left-hand square and the whole rectangle are pullbacks, the right-hand square is also a pullback. A regular category $\Ac$ is called \emph{(Barr) exact} when any equivalence relation in $\Ac$ is effective (i.e.\ a kernel pair) \cite{Barr}. A category $\Ac$ is {\em pointed} when it has a zero object $0$ (i.e.\ an initial object that is also terminal).
 
\begin{definition}\cite{Janelidze-Marki-Tholen}
A pointed category $\Ac$ is {\em semi-abelian} when it is exact, it has binary coproducts and it is protomodular.
\end{definition}

A characterization of the algebraic theories with the property that the corresponding category of algebras is a semi-abelian category was obtained by Bourn and Janelidze:

\begin{theorem}\cite{Bourn-Janelidze}
 A variety of universal algebras $\mathcal V$ is semi-abelian if and only if its theory $\mathbb{T}$ has a unique constant $0$, binary terms $t_1, \dots , t_n$ and a $(n+1)$-ary term $\tau$ satisfying the identities $\tau (x,t_1 (x,y), \dots, t_n (x,y))=y$ and $t_i(x,x)=0$ for each $i=1,\dots, n$.\noproof 
\end{theorem}

Classical examples of semi-abelian varieties are groups, rings, Lie algebras, commutative algebras, crossed modules, precrossed modules and Heyting semilattices \cite{Jo}. Compact Hausdorff (profinite) groups or, more generally, compact Hausdorff (profinite) semi-abelian algebras are semi-abelian categories \cite{Borceux-Clementino}, as is the dual of the category of pointed sets, or the category of $C^{\star}$-algebras.
 
An important property of semi-abelian categories is the fact that every regular epimorphism $f \colon  X \to Y$ is the cokernel of its kernel \cite{Bourn2001}. In other words, the class of regular epimorphisms coincides with the class of normal epimorphisms. In a semi-abelian category, a short exact sequence 
\[
\xymatrix{0 \ar[r] & K \ar@{{ >}->}[r]^-{k} & X  \ar@{-{ >>}}[r]^-{f}   &  Y \ar[r] & 0}
\]
is then a zero sequence ($f \comp k =0$, with $0$ the zero arrow) such that $f$ is a regular epi and $k$ is a kernel of $f$.
Semi-abelian categories are known to provide an appropriate setting for the description of some important aspects of homological algebra, modelled on the category of groups \cite{Borceux-Bourn}. 

The next property, due to Bourn, will be needed in what follows:

\begin{proposition}\cite{Bourn1991, Bourn2001} \label{PropositionLeftRightPullbacks}
In a semi-abelian category $\Ac$, let us consider the commutative diagram with exact rows 
\begin{equation}\label{Diagram-Exact-Rows}
\vcenter{\xymatrix{0 \ar[r] & K' \ar[d]_{u} \ar@{{ >}->}[r]^{k'} & X' \ar[d]^{v} \ar@{-{ >>}}[r]^{f'}   &  Y'\ar[d]^{w} \ar[r] & 0\\
0 \ar[r] & K \ar@{{ >}->}[r]_{k} & X  \ar@{-{ >>}}[r]_{f}   &  Y \ar[r] & 0. }}
\end{equation}
Then:
\begin{enumerate}
\item $u$ is an isomorphism if and only if the right-hand square is a pullback;
\item $w$ is a monomorphism if and only if the left-hand square is a pullback.
\end{enumerate}
\end{proposition}
\begin{proof}
1. Given the commutative diagram
\[
\xymatrix{ K' \ar[d]_{}  \ar@{}[dr]|-{\texttt{(i)}}\ar[r]^{k'} & X' \ar@{-{ >>}}[d]^{f'} \ar@{}[dr]|-{\texttt{(ii)}} \ar[r]^{v}   &  X \ar[d]^{f}  \\
0 \ar[r]_{} & Y'  \ar[r]_{w}   &  Y, }
\]
one has that \texttt{(i)} is a pullback by construction. But $v\comp k'=k \comp u$, so that the whole rectangle \texttt{(i)}+\texttt{(ii)} is a pullback whenever $u$ is an isomorphism. From the fact that $\Ac$ is semi-abelian and $f'$ is a regular epimorphism it follows that \texttt{(ii)} is a pullback.

Conversely, if we assume that the square \texttt{(ii)} is a pullback, so is the following rectangle:
\[
\xymatrix{ K' \ar[d]_{} \ar[r]^{u} & K \ar[d]^{}  \ar[r]^{k}   &  X \ar[d]^{f}  \\
0 \ar@{=}[r]_{} & 0  \ar[r]   &  Y. }
\]
Since the right-hand square is a pullback by construction, it follows that the left-hand square is a pullback as well, and $u$ is an isomorphism.

2. The non-trivial implication essentially follows from the fact that, in a semi-abelian category, an arrow is a monomorphism if its kernel is $0$ \cite{Bourn1991}.
\end{proof}

The following property will be also needed:

\begin{proposition} \label{Proposition-pushout}
In a semi-abelian category $\Ac$, let us consider the commutative diagram with exact rows (\ref{Diagram-Exact-Rows}). If the left-hand square is a pushout, then $w$ is an isomorphism; conversely, if $w$ is an isomorphism and $u$ is a regular epimorphism, then the left-hand square is a pushout.
\end{proposition}
\begin{proof}
If the left-hand side square is a pushout then $f'$ and $0\colon K\to Y'$ induce an arrow $f''\colon X\to Y'$ satisfying $f''\comp v=f'$ and $f''\comp k=0$. Then also $w\comp f''=f$. Moreover, like $f$, $f''$ is a cokernel of $k$; hence the unique comparison map $w$ is an isomorphism.

Now again consider the diagram (\ref{Diagram-Exact-Rows}); pushing out $u$ along $k'$, then taking a cokernel $f''$ of $\overline{k'}$ induces the dotted arrows in the diagram below.
\[
\xymatrix{0 \ar[r] & K' \ar@{-{ >>}}[d]_{u} \ar@{{ >}->}[r]^{k'} \ar@{}[rd]|>{\pushout}& X' \ar@{.{ >>}}[d]^{\overline{u}} \ar@{-{ >>}}[r]^{f'}   &  Y'\ar@{:}[d] \ar[r] & 0\\
0 \ar@{.>}[r] & K \ar@{:}[d] \ar@{.>}[r]_-{\overline{k'}} & P \ar@{.>}[d]^{v'} \ar@{.{ >>}}[r]_-{f''}   &  Y'\ar[d]^{w} \ar@{.>}[r] & 0\\
0 \ar[r] & K \ar@{{ >}->}[r]_{k} & X  \ar@{-{ >>}}[r]_{f}   &  Y \ar[r] & 0}
\]
Note that $\overline{k'}$ is a monomorphism because $k$ is one; being the regular image of the kernel $k'$, $\overline{k'}$ is a kernel as well~\cite{Bourn2003}. The Short Five Lemma implies that $w$ is an isomorphism if and only if so is $v'$.
\end{proof}

We conclude this section by recalling the following definition:

\begin{definition} \cite{Carboni-Lambek-Pedicchio}
A finitely complete category $\Ac$ is a {\em Mal'tsev category} if every internal reflexive relation in $\Ac$ is an equivalence relation.
\end{definition}

Thanks to a result of Bourn, it is well-known that every protomodular category is a Mal'tsev category.

\section{Centrality}\label{Section-Centrality}

In this section we explore different definitions of centrality. The first notion is the classical notion of {\em centrality of congruences} introduced by Smith in the context of Mal'tsev varieties \cite{Smith}, which has later been extended to Mal'tsev categories \cite{Carboni-Pedicchio-Pirovano, Pedicchio}. The second one is the notion of  \emph{central arrows}, first defined by Huq \cite{Huq} in a context that is essentially equivalent to the context of semi-abelian categories \cite{Janelidze-Marki-Tholen}. In any pointed protomodular category, there is a natural way to compare the two notions, because in such a category normal subobjects correspond to (internal) equivalence relations. The paper \cite{BG} was the first in which the relationship between these two notions of centrality was investigated. Here we present two new results in this direction, Proposition~\ref{Proposition-Huq-Centrality-is-Smith-Centrality} and Proposition~\ref{Proposition-Subobject-of-Central-is-Normal}, which will also be useful in the subsequent sections.

Let us begin by recalling the following definition due to Bourn.

\begin{definition}\label{Definition-Normal-Monomorphism} \cite{Bourn2000}
An arrow $k\colon  K \to X$ in a finitely complete category $\Ac$ is \textit{normal to an equivalence relation} $R$ on $X$ when: (1) $k^{-1}(R)$ is the largest equivalence relation $\nabla_K$ on $K$; (2) the induced map $\nabla_K \to R$ in the category $\Eq (\Ac)$ of internal equivalence relations in $\Ac$ is a discrete fibration.
\end{definition}

This means that
\begin{enumerate}
\item there is a map $\tilde k\colon K \times K \to R$ in $\Ac$ such that the diagram
\[
\xymatrix{ K \times K \ar@{.>}[r]^-{\tilde{k}} \ar@<.8ex>[d]^-{\pi_{2}} \ar@<-.8ex>[d]_-{\pi_{1}} & R \ar@<.8ex>[d]^-{\pi_{2}} \ar@<-.8ex>[d]_-{\pi_{1}}\\
K \ar[r]_-{k} & X}
\]
commutes;
\item any of the commutative squares in the diagram above is a pullback.
\end{enumerate}

It can be proved that the arrow $k$ is then necessarily a monomorphism; furthermore, when the category $\Ac$ is protomodular, a monomorphism can be normal to at most one equivalence relation, so that the fact of being normal becomes a property \cite{Bourn2000}. The notion of normal monomorphism gives an intrinsic way to express the fact that $K$ is an equivalence class of $R$.

In a pointed finitely complete category $\Ac$ there is a natural way to associate, with any equivalence relation 
\[
\xymatrix{{R}  \ar@<-1.6ex>[rr]_-{\pi_{1}} \ar@<1.6ex>[rr]^-{\pi_{2}} && {X}, \ar[ll]|-{\Delta}}
\]
a normal subobject $k_{R}$, called the \emph{normalization} of $R$, or \emph{the normal subobject associated with $R$}: it is defined as the composite $k_{R}=\pi_{2}\comp \ker \pi_{1}$
\[
K[\pi_{1}] \to/{ >}->/^{\ker \pi_{1}} R \to^{\pi_{2}} X.
\]
In the pointed protomodular case, this construction determines a bijection between the equivalence relations on $X$ and the normal subobjects of $X$ \cite{Bourn2000}.

Two equivalence relations $R$ and $S$ on an object $X$ \emph{centralize (in the sense of Smith)} when there exists a double equivalence relation $C$ on $R$ and $S$ such that any commutative square in the diagram
\[
\xymatrix{C \ar@<-.8ex>[d]_{p_1} \ar@<.8ex>[d]^{p_2} \ar@<-.8ex>[r]_{p_2} \ar@<.8ex>[r]^{p_1} & S \ar@<-.8ex>[d]_-{\pi_{1}} \ar@<.8ex>[d]^-{\pi_{2}}\\
R \ar@<-.8ex>[r]_-{\pi_{2}} \ar@<.8ex>[r]^-{\pi_{1}} & X}
\]
is a pullback \cite{Smith, Carboni-Pedicchio-Pirovano}. In this case, $C$ is called the centralizing double relation on $R$ and $S$. An equivalence relation $R$ on $X$ is said to be \emph{central} when $R$ and $\nabla_{X}$ (the largest equivalence relation on $X$) centralize.

Adopting the terminology due to Bourn \cite{Bourn2002}, we say that two coterminal morphisms $k\colon K\to X$ and $k'\colon K'\to X$ in a pointed finitely complete category \emph{cooperate} when a morphism $\varphi_{k,k'}\colon K\times K'\to X$ exists satisfying $\varphi_{k,k'}\comp l_{K}=k$ and $\varphi_{k,k'}\comp r_{K'}=k'$, where $l_K=(1_K,0) \colon  K \to K \times K'$ and $r_{K'} = (0,1_{K'}) \colon  K' \to K \times K'$. (In Huq's terminology, $k$ and $k'$ \emph{commute} \cite{Huq}.) The arrow $\varphi_{k,k'}$ is called a \emph{cooperator} of $k$ and $k'$. In particular, an arrow $k$ is said to be \emph{central (in the sense of Huq)} when $k$ and $1_{X}$ cooperate.
 
It is well-known that, in general, two equivalence relations $R$ and $S$ need not centralize when $k_{R}$ and $k_{S}$ cooperate, not even in a variety of $\Omega$-groups (see \cite{Bourn2004} for a counter-example). However, we are now going to show that this \emph{is} the case in any pointed protomodular category, whenever $R$ (or $S$) is $\nabla_{X}$:

\begin{proposition}\label{Proposition-Huq-Centrality-is-Smith-Centrality}
Let $\Ac$ be a pointed protomodular category. 
An equivalence relation $R$ in $\Ac$ is central if and only if its associated normal subobject $k_{R}$ is central.
\end{proposition}
\begin{proof}
Let us first assume that the equivalence relation 
\[
\xymatrix{{R}  \ar@<-1.6ex>[rr]_-{\pi_{1}} \ar@<1.6ex>[rr]^-{\pi_{2}} && {X} \ar[ll]|-{\Delta}}
\]
is central, and let $C$ be the associated centralizing double relation on $R$ and $\nabla_{X}$. We can then consider the following diagram
\[
\xymatrix{C \ar[d]_{p_1}  \ar@<-.8ex>[r]_{p_2} \ar@<.8ex>[r]^{p_1} & R \ar@{}[rd]|-{\texttt{(i)}} \ar[d]_-{\pi_{1}} \ar@{.>>}[r] & K  \ar@{.>}[d] \\
X \times X \ar@<-.8ex>[r]_(.6){\pi_{2}} \ar@<.8ex>[r]^(.6){\pi_{1}}  & X \ar@{.>>}[r] & 0,}
\]
where the square \texttt{(i)} is obtained by taking coequalizers. Since $C$ is centralizing, both left hand side squares are pullbacks, hence so is \texttt{(i)}, and $R \cong X \times K$. This moreover induces the commutative diagram
\[
\xymatrix{K \ar[d] \ar[r]^-{k} & R \ar@{}[rd]|<{\pullback} \ar[d]_-{\pi_{1}} \ar[r]^-{p} & K \ar[d] \\
0 \ar[r]  & X \ar[r] & 0}
\]
where $p\comp k=1_{K}$, and both the outer rectangle and the right hand square are pullbacks. It follows that $K$ is the normal subobject associated with $R$. By considering also the second projection $\pi_2$ from $R$ to $X$ one can easily check that $R$ is then canonically isomorphic to the equivalence relation
\begin{equation}\label{Diagram-Centralizing-Relation}
\xymatrix{{K\times X} \ar@<-1.6ex>[rr]_-{\pi_{X}} \ar@<1.6ex>[rr]^-{\varphi_{k_{R},1_{X}}} && {X}; \ar[ll]|-{r_{X}}}
\end{equation}
the arrow $\varphi_{k_{R},1_{X}}$ is the needed cooperator.

Conversely, let us suppose that $k$ is a central monomorphism with cooperator $\varphi_{k,1_{X}}$. One can then form the reflexive graph (\ref{Diagram-Centralizing-Relation})---call it $R_{k}$. It is a relation, since the commutative square
\[
\xymatrix{K \ar[d]_{k}  \ar[r]^-{l_K} & K \times X \ar[d]^{(\varphi_{k,1_{X}}, \pi_2)} \\
X \ar[r]_-{l_{X}}  & X \times X}
\]
is a pullback, and in a protomodular category pullbacks reflect monomorphisms \cite{Bourn1991}.
Since $R_{k}$ is a reflexive relation in a Mal'tsev category, it is an equivalence relation. Furthermore, this equivalence relation corresponds to $k$ via the bijection between normal subobjects and equivalence relations: to see this, it suffices to observe that  the normalization of $R_{k}$ is isomorphic to the normalization of the opposite relation $R_{k}^{\op}$.

Finally, consider the double equivalence relation determined by the kernel pair $R[\pi_{K}]$ of $\pi_{K}\colon K\times X\to K$:
\[
\xymatrix{R[\pi_{K}] \ar@<-.8ex>[d] \ar@<.8ex>[d] \ar@<-.8ex>[r] \ar@<.8ex>[r] & X\times X \ar@<-.8ex>[d]_{\pi_1} \ar@<.8ex>[d]^{\pi_2}\\
K\times X \ar@<-.8ex>[r]_{\pi_{X}} \ar@<.8ex>[r]^{\varphi_{k,1_{X}}} & X.}
\]
It is clearly a 
centralizing double relation on $R_{k}$ and $\nabla_{X}$, as desired.
\end{proof}

Recall that an \emph{extension of an object $Y$ (by an object $K$)} is a regular epimorphism $f\colon X\to Y$ with its kernel $K$:
\[
0 \to K \to/{{ >}->}/^{k} X \to/-{ >>}/^{f} Y \to 0.
\]
The category of extensions of $Y$ (considered as a full subcategory of the slice category $(\Ac\downarrow Y)$) is denoted by $\Ext (Y)$; the category of all extensions in $\Ac$ (considered as a full subcategory of the arrow category $\Fun (\mathsf{2},\Ac)$: morphisms are commutative squares) by $\Ext (\Ac)$. Recall that in a semi-abelian category, a subobject is normal if and only if it is a kernel. An extension $f\colon X\to Y$ is called \emph{central} if its kernel is central in the sense of Huq, i.e.\ if $\ker f$ cooperates with $1_{X}$. We write $\CentrExt (Y)$ for the full subcategory of $\Ext (Y)$ determined by the central extensions. The following well-known property of central extensions of groups will now be shown to hold in any semi-abelian category:

\begin{proposition}\label{Proposition-Subobject-of-Central-is-Normal}
Let $\Ac$ be a semi-abelian category and $f\colon X\to Y$ a central extension in $\Ac$. Every subobject of $k=\ker f\colon K\to X$ is normal in $X$. 
\end{proposition}
\begin{proof}
Let $i\colon M\to K$ be a monomorphism and denote $m=k\comp i$. We are to show that $m$ is a normal monomorphism. Now $k$ is central, hence so is $m$. This means that $m$ cooperates with $1_{X}$; in particular, there exists an arrow $\varphi_{m,1_{X}}\colon M\times X\to X$ satisfying $\varphi_{m,1_{X}}\comp l_{M}=m$ (where $l_{M}=(1_M,0)$). But the arrow $l_{M}$, being the kernel of $\pi_{2}$, is normal; hence so is $m$, since it is the regular image of $\varphi_{m,1_{X}}\comp l_{M}$ \cite{Bourn2003}.
\end{proof}

\section{Abelian objects and central extensions}\label{Section-Abelianization}

In this section, $\Ac$ will be a semi-abelian category. An object $X$ in $\Ac$ is called \emph{abelian} when there is a centralizing relation on $\nabla_{X}$ and $\nabla_{X}$. The full subcategory of $\Ac$ determined by the abelian objects is denoted $\Ab(\Ac)$. Since $\Ac$ is pointed, $\Ab(\Ac)$ coincides with the category of internal abelian group objects in $\Ac$. $\Ab(\Ac)$ is an abelian category. It is well-known that $\Ab(\Ac)$ is a reflective subcategory of $\Ac$, closed in $\Ac$ under subobjects and quotients (i.e.\ it is a Birkhoff subcategory) \cite{Bourn-Gran}.
\[
\xymatrix{\Ac \ar@<.8ex>[r]^-{\ab} \ar@{}[r]|-{\perp} & {\,\Ab(\Ac)} \ar@<.8ex>@{_(->}[l]}
\]
For any object $X$, we shall denote the $X$-component of the unit of the adjunction by $\eta_{X}\colon X\to X/[X,X]=\ab (X)$, and its kernel by $\mu_{X}\colon [X,X]\to X$, the $X$-component of a natural transformation $\mu \colon V\To 1_{\Ac }\colon \Ac \to \Ac$.

Similarly, for any object $Y$, the category $\CentrExt (Y)$ is reflective and closed under subobjects and regular quotients in $\Ext (Y)$.
\[
\xymatrix{\Ext (Y) \ar@<.8ex>[r]^-{\centr} \ar@{}[r]|-{\perp} & {\,\CentrExt (Y)} \ar@<.8ex>@{_(->}[l]}
\]
The $f$-component of the unit of the adjunction is given by the horizontal arrows in the diagram
\[
\xymatrix{X \ar@{-{ >>}}[r] \ar@{-{ >>}}[d]_-{f} & \textstyle{\frac{X}{[K,X]}} \ar@{-{ >>}}[d]^-{\centr f}\\
Y \ar@{=}[r] & Y,}
\]
where $k\colon K\to X$ denotes a kernel of $f\colon X\to Y$ \cite[Theorem 2.8.11]{Borceux-Bourn}.

\begin{proposition}\label{Proposition-Split-Mono-Central-Extension}
Consider the diagram of solid arrows
\[
\xymatrix{0 \ar[r] & K \ar@{{ >}->}[r]^-{k} \ar@{{ >}->}@<-.8ex>[d]_-{m} & X \ar@{-{ >>}}[r]^-{f} \ar@{{ >}.>}@<-.8ex>[d]_-{\overline{m}} & Y \ar@{:}[d] \ar[r] & 0\\
0 \ar@{.>}[r] & A \ar@{-{ >>}}@<-.8ex>[u]_-{s} \ar@{{ >}.>}[r] & Z \ar@{.>>}@<-.8ex>[u] \ar@{.>>}[r] & Y \ar@{.>}[r] & 0.}
\]
If the above sequence is exact, $m$ is a normal monomorphism split by $s$, $f$ is a central extension and $A$ is abelian, then a central extension of $Y$ by $A$ exists making the diagram commutative.
\end{proposition}
\begin{proof}
Let $q\colon A\to Q$ denote a cokernel of $m$, and let us consider the sequence
\[
\xymatrix{0 \ar[r] & K \ar@<-.8ex>@{{ >}->}[r]_-{m} & A \ar@<-.8ex>@{-{ >>}}[l]_-{s} \ar@{-{ >>}}[r]^-{q} & Q \ar[r] & 0}
\]
in $\Ab(\Ac)$, which is a split exact sequence.  It follows that $A$ is a product of $K$ with $Q$ and that, up to isomorphism, $m$ is $l_{K}\colon K\to K\times Q$ and $s$ is $\pi_{1}\colon {K\times Q\to K}$. In the diagram
\[
\xymatrix{0 \ar[r] & K \ar@{{ >}->}[r]^-{k} \ar@{}[rd]|-{\texttt{(i)}} \ar@{{ >}->}@<-.8ex>[d]_-{l_{K}} & X \ar@{-{ >>}}[r]^-{f} \ar@{{ >}->}@<-.8ex>[d]_-{l_{X}} & Y \ar@{=}[d] \ar[r] & 0\\
& K\times Q \ar@{-{ >>}}@<-.8ex>[u]_-{\pi_{1}} \ar@{{ >}->}[r]_-{k\times 1_{Q}} & {X\times Q} \ar@{-{ >>}}@<-.8ex>[u]_-{\pi_{1}} \ar@{-{ >>}}[r]_-{f\comp \pi_{1}} & Y}
\]
the upward-pointing square \texttt{(i)} is a pullback; it follows that $k\times 1_{Q}$ is a kernel of $f\comp \pi_{1}$, and then $f\comp \pi_{1}$ is the cokernel of $k\times 1_{Q}$. From the fact that $Q$ is abelian and $k$ is central one concludes that $k\times 1_{Q}$ is central, and this completes the proof.
\end{proof}

\begin{proposition}\label{Proposition-Reg-Epi-Central-Extension}
Consider the diagram of solid arrows
\[
\xymatrix{0 \ar[r] & K \ar@{{ >}->}[r]^-{k} \ar@{-{ >>}}[d]_-{r} & X \ar@{-{ >>}}[r]^-{f} \ar@{.>>}[d]^-{\overline{r}} & Y \ar@{:}[d] \ar[r] & 0\\
0 \ar@{.>}[r] & A \ar@{{ >}.>}[r] & Z \ar@{.>>}[r] & Y \ar@{.>}[r] & 0.}
\]
If the above sequence is exact, $r$ is a regular epimorphism, $f$ is a central extension and $A$ is abelian, then a central extension of $Y$ by $A$ exists making the diagram commutative.
\end{proposition}
\begin{proof}
We may define $\overline{r}$ as a cokernel of $k\comp \ker r$:
\[
\xymatrix{& K[r] \ar@{=}[r] \ar@{{ >}->}[d]_-{\ker r} & K[r] \ar@{{ >}->}[d]^-{k\comp \ker r} \\
0 \ar[r] & K \ar@{}[rd]|{\texttt{(i)}} \ar@{{ >}->}[r]^-{k} \ar@{-{ >>}}[d]_-{r} & X \ar@{.{ >>}}[d]_-{\overline{r}} \ar@{-{ >>}}[r]^-{f} & Y \ar[r] \ar@{:}[d] & 0\\
0 \ar@{.>}[r] & A \ar@{{ >}.>}[r]_-{\overline{k}} & Z \ar@{.{ >>}}[r] & Y \ar@{.>}[r] & 0.}
\]
The morphism $k\comp \ker r$ is a kernel thanks to Proposition~\ref{Proposition-Subobject-of-Central-is-Normal}. Taking cokernels induces the square \texttt{(i)}, which is easily seen to be a pushout. Remark also that, by Proposition~\ref{PropositionLeftRightPullbacks}, $\overline{k}$ is a monomorphism, hence a kernel, since it is the regular image of $k$ along $\overline{r}$ \cite{Bourn2003}. Taking a cokernel of $\overline{k}$ gives rise to the rest of the diagram, thanks to Proposition~\ref{Proposition-pushout}. The induced extension is central, because $\CentrExt (\Ac)$ is closed under quotients in $\Ext (\Ac)$---and a quotient in $\Ext (\Ac)$ is a pushout in $\Ac$ of a regular epimorphism along a regular epimorphism.
\end{proof}

\begin{corollary}\label{Corollary-Central-Extension-Pushout}
Consider the diagram of solid arrows
\[
\xymatrix{0 \ar[r] & K \ar@{{ >}->}[r]^-{k} \ar[d]_-{a} & X \ar@{-{ >>}}[r]^-{f} \ar@{.>}[d]^-{\overline{a}} & Y \ar@{:}[d] \ar[r] & 0\\
0 \ar@{.>}[r] & A \ar@{{ >}.>}[r]_-{k_{a}} & Z_{a} \ar@{.>>}[r]_-{p_{a}} & Y \ar@{.>}[r] & 0.}
\]
If the above sequence is exact, $f$ is a central extension and $A$ is abelian, then a central extension of $Y$ by $A$ exists making the diagram commutative.
\end{corollary}
\begin{proof}
Since $K$ and $A$ are abelian, the arrow $a\colon K\to A$ lives in the abelian category $\Ab(\Ac)$, and thus may be factored as the composite
\[
K\to/{ >}->/^{l_{K}} K\oplus A \to/->>/^{[a,1_{A}]} A. 
\]
The result now follows from the previous propositions. 
\end{proof}

Alternatively, this result may be obtained by using torsor theory \cite{Bourn-Janelidze:Torsors}.

\section{The perfect case: universal central extensions}\label{Section-Perfect-Case}

Suppose that $\Ac$ is a semi-abelian category with enough (regular) projectives and $Y$ an object of $\Ac$. Then the category $\CentrExt (Y)$ always has a weakly initial object: for if $f\colon X\to Y$ is a \emph{(projective) presentation} of $Y$, i.e.\ a regular epimorphism with $X$ projective, then the reflection $\centr f\colon X/[K,X]\to Y$ of $f$ into $\CentrExt (Y)$ is a central extension of $Y$. It is weakly initial, as any other central extension $g\colon Z\to Y$ induces a morphism $\centr f\to g$ in $\CentrExt (Y)$, the object $X$ being projective.

An initial object in $\CentrExt (Y)$ is called a \emph{universal} central extension of $Y$. In contrast with the existence of weakly initial objects, for a universal central extension of $Y$ to exist, the object $Y$ must be \emph{perfect}: such is an object $Y$ of $\Ac$ with the property that its reflection $\ab (Y)$ into $\Ab(\Ac)$ is $0$, i.e. $[Y,Y]\cong Y$.

To see this, let us first recall the Fundamental Theorem of Categorical Galois Theory in a form which is suitable for our context (see \cite{Janelidze-Kelly}). Consider a projective presentation $f \colon {X \to Y}$. Its \emph{Galois groupoid} ${\Gal(f)}$ is the image, under the functor $\ab \colon {\Ac \to \Ab(\Ac)}$, of the equivalence relation 
$$\xymatrix@=30pt{ 
 R[f] \times_X R[f] \ar@<3ex>[r]^-{\pi_{1}}
\ar@<-3ex>[r]^-{\pi_{2}} 
  \ar@<0ex>[r]^-{m} &  {R[f]\, } \ar@<3ex>[r]^-{p_{1}} \ar@<-3ex>[r]^-{p_{2}} &  X \ar@<0ex>[l]_-{\Delta }}$$
that is the kernel relation of $f$. Remark that, in our situation, the diagram ${\Gal(f)}$ is indeed an internal groupoid in $\Ab(\Ac)$: as explained in \cite{Gran-Alg-Cent}, this is a consequence of the fact that the reflector of a semi-abelian category to a Birkhoff subcategory preserves any pullback of a split epimorphism along a split epimorphism.

By taking into account the main result of \cite{Bourn-Gran}, the Fundamental Theorem of Galois Theory may be written as a category equivalence 
\[
\CentrExt(Y) \simeq \{{\Gal(f)}, \Ab(\Ac) \}
\]
between the category $\CentrExt (Y)$ of central extension of $Y$ and the category $\{{\Gal(f)}, {\Ab(\Ac)} \}$ of discrete fibrations in $\Ab(\Ac)$ on the Galois groupoid ${\Gal(f)}$, whose objects can be represented by diagrams of the form 
\begin{equation}\label{discrete}
\vcenter{ \xymatrix@=50pt{G_{2} \ar@<3ex>[r]^-{\pi_{1}} \ar[d]^-{g_2}
\ar@<-3ex>[r]^-{\pi_{2}} 
  \ar@<0ex>[r]^-{m} &
 G_{1} \ar[d]^{g_1} \ar@<3ex>[r]^{d_{1}} \ar@<-3ex>[r]^{d_{2}} & G_{0} \ar[d]^{g_0} \ar@<0ex>[l]_{e} \\
 \ab(R[f] \times_X R[f]) \ar@<3ex>[r]^-{\ab(\pi_{1})}
\ar@<-3ex>[r]^-{\ab(\pi_{2})} 
  \ar@<0ex>[r]^-{\ab(m)} &
 \ab(R[f]) \ar@<3ex>[r]^{\ab( p_{1})} \ar@<-3ex>[r]^{\ab(p_{2})} & \ab (X) \ar@<0ex>[l]_{\ab(\Delta )}}}
\end{equation}
with the property that the arrow $g_0 \colon G_0 \to \ab (X)$ is a regular epimorphism.
We are now ready to prove the following

\begin{proposition}\label{Proposition-Universal-Central-Extension}
Let $\Ac$ be a semi-abelian category with enough projectives. An object $Y$ of $\Ac$ is perfect if and only if $Y$ admits a universal central extension.
\end{proposition}
\begin{proof}
First suppose that $\CentrExt (Y)$ has an initial object $u\colon U\to Y$. Because $\pi_{1}\colon Y\times \ab (Y)\to Y$ is central, a unique morphism $(u,y)\colon  U\to Y\times \ab (Y)$ exists. But then $0\colon U\to \ab (Y)$ is equal to $\eta_{Y}\comp u\colon U\to \ab (Y)$, and $\ab (Y)=0$.

Conversely, consider a presentation $f \colon X\to Y$ of a perfect object $Y$. Its Galois groupoid $\Gal(f)$ is connected, since the functor $\ab \colon \Ac \to \Ab(\Ac)$ preserves coequalizers, as a left adjoint.

Now, in the abelian category $\Ab(\Ac)$, the normalization functor recalled in Section~\ref{Section-Centrality} determines a category equivalence between the internal groupoids in $\Ac$ and the arrows of $\Ac$, which restricts to a category equivalence between the internal connected groupoids in $\Ac$ and the epimorphisms. Let $\phi \colon {K_f \to \ab(X)}$ be the ``normalization'' of $\Gal(f)$: by applying Proposition~\ref{PropositionLeftRightPullbacks}.1.\ to diagram (\ref{discrete}) one can check that the category $\{{\Gal(f)}, {\Ab(\Ac)} \}$ is equivalent to the category of triples $(Z,k,\pi)$, where $k \colon {K_f \to Z}$ is a monomorphism and $\pi \colon {Z \to \ab(X)}$ is an epimorphism. It follows that, when $\Gal(f)$ is connected, so that its normalization $\phi \colon {K_f \to \ab(X)}$ is an epimorphism, this category has an initial object, namely $(K_f, 1_{K_f}, \phi)$.
\end{proof}

\section{Cohomology}\label{Section-Cohomology}

 From now on, $\Ac$ will be a semi-abelian category, monadic over the category $\Set$ of sets. Such categories were characterized by Gran and Rosick\'y in \cite{Gran-Rosicky:Monadic} by extending a previous result due to Bourn and Janelidze \cite{Bourn-Janelidze}. We recall some concepts, results and notation from the paper \cite{EverVdL2}, adapted to our present situation.

Let  
\[
\Gb =(G\colon \Ac \to \Ac ,\quad\epsilon\colon G\To 1_{\Ac} ,\quad \delta\colon G\To G^{2})
\]
denote the comonad on $\Ac$, induced by the monadicity requirement. Recall that the axioms of comonad state that, for every object $X$ of $\Ac$, $\epsilon_{GX}\comp \delta_{X}=G\epsilon_{X}\comp \delta_{X}=1_{GX}$ and $\delta_{GX}\comp \delta_{X}=G\delta_{X}\comp \delta_{X}$. Putting
\[
\del_{i}=G^{i}\epsilon_{G^{n-i}X}\colon G^{n+1}X\to G^nX
\]
and
\[
\sigma_{i}=G^{i}\delta_{G^{n-i}X}\colon G^{n+1}X\to G^{n+2}X ,
\]
for $0\leq i \leq n$, makes the sequence $(G^{n+1}X)_{n\in \N}$ a simplicial object $\Gb X$ of $\Ac$. This induces a functor from $\Ac$ to the category $\simpA$ of simplicial objects in $\Ac$: this functor will also be denoted $\Gb$.

\begin{remark}\label{Remark-Projectives}
A consequence of the monadicity of $\Ac$ is the existence of sufficiently many regular projective objects in $\Ac$. Indeed, any value $G (X)$ of $G$ is projective, and the morphism $\epsilon_{X}\colon G (X)\to X$ is a regular epimorphism. By calling \emph{(projective) presentation of} $Y$ an extension $f\colon X\to Y$ with $X$ projective, it follows that, for every object $Y$, such a projective presentation exists.
\end{remark}

Recall that a chain complex in a semi-abelian category is called \emph{proper} when its differentials have normal images. As in the abelian case, the $n$-th homology object $H_{n}C$ of a proper chain complex $C$ with differentials $d_n$ is the cokernel of $C_{n+1}\to K[d_n]$. The \emph{normalization functor} $N\colon \simpA \to \PCh \Ac$ turns a simplicial object $A$ into the \emph{Moore complex} $N (A)$ of $A$, the chain complex with $N_{0}A=A_{0}$,
\[
N_{n} A=\bigcap_{i=0}^{n-1}K[\del_{i}\colon A_{n}\to A_{n-1}]
\]
and differentials $d_{n}=\del_{n}\comp \bigcap_{i}\ker \del_{i}\colon N_{n} A\to N_{n-1} A$, for $n\geq 1$, and $A_{n}=0$, for $n<0$. Since $N (A)$ is a proper chain complex in a semi-abelian category, one can define its homology objects in the usual way.

\begin{definition}\cite[Section 6, case $\Bc =\Ab(\Ac)$]{EverVdL2}\label{Definition-Homology}
For $n\in \N_{0}$, the object
\[
H_{n} X=H_{n-1} N \ab (\Gb X)
\]
is the \emph{$n$-th homology object of $X$ (with coefficients in $\ab $) relative to the cotriple $\Gb$}. This defines a functor $H_{n} \colon \Ac \to \Ab(\Ac) $, for any $n\in \N_{0}$. 
\end{definition}

\begin{proposition}\label{Proposition-H_0-Sequence}
Let
\[
\xymatrix{ 0 \ar[r] & K \ar@{{ >}->}[r] & X \ar@{-{ >>}}[r]^-{f} & Y \ar[r] & 0 }
\]
be a short exact sequence in $\Ac$. Then the induced sequence in $\Ab (\Sc \Ac)$ 
\begin{equation}\label{Sequence-simpB}
 0 \to K[\ab\Gb f] \to/{{ >}->}/^{\ker \ab\Gb f} \ab\Gb X \to/{-{ >>}}/^{\ab\Gb f} \ab\Gb Y \to 0
\end{equation}
is degreewise split exact and is such that 
\[
H_{0}K[\ab \Gb f]\cong\textstyle{\frac{K}{[K,X]}}.
\]
\end{proposition}
\begin{proof}
Since $G$ turns regular epimorphisms into split epimorphisms, the simplicial morphism $\ab\Gb f$ is degreewise split epimorphic in $\Ab (\Sc \Ac)$.

For any $n\geq 1$, the short exact sequence 
\[
0 \to K[G^{n} f] \to/{{ >}->}/^{\ker G^{n} f} G^{n} X \to/{-{ >>}}/^{G^{n} f} G^{n} Y \to 0,
\]
through \cite[Theorem 5.9]{EverVdL1}, induces the exact sequence
\[
0 \to \tfrac{K[G^{n} f]}{[K[G^{n} f], G^{n}X]} \to/{{ >}->}/^{\ker \ab G^{n} f} \ab G^{n} X \to/{-{ >>}}/^{\ab G^{n} f} \ab G^{n} Y \to 0.
\]
As a consequence, 
\[
K[\ab \Gb f]\cong \tfrac{K[\Gb f]}{[K[\Gb f], \Gb X]}. 
\]
It is now possible to prove $H_{0}K[\ab \Gb f]\cong K/[K,X]$ by showing that the fork
\[
\xymatrix{{\tfrac{K[G^{2}f]}{[K[G^{2}f],G^{2}X]}} \ar@<-.8ex>[r] \ar@<.8ex>[r] & {\tfrac{K[G f]}{[K[G f],G X]}} \ar[r] & {\tfrac{K}{[K,X]}}}
\]
is a coequalizer diagram \cite{VdLinden:Doc}.
\end{proof}

\begin{theorem}[Stallings-Stammbach sequence and Hopf formula]\label{Theorem-Homology} \cite{EverVdL2}
If
\[
\xymatrix{ 0 \ar[r] & K \ar@{{ >}->}[r] & X \ar@{-{ >>}}[r]^-{f} & Y \ar[r] & 0 }
\]
is a short exact sequence in $\Ac$, then there exists an exact sequence
\begin{equation}\label{Stallings-Stammbach-Sequence}
 H_2 X \to<250>^{H_2 f} H_{2} Y \to<250> \textstyle{\frac{K}{[K,X]}} \to<250> H_{1} X \to/{-{ >>}}/<250>^{H_{1}f} H_{1} Y \to<250> 0
\end{equation}
in $\Ab(\Ac)$ that depends naturally on the given short exact sequence. Moreover $H_{1} Y\cong \ab (Y)$ and, when $X$ is projective, $H_{2} Y\cong {(K\cap [X,X])}/{[K,X]}$.\noproof 
\end{theorem}

\begin{remark}\label{Remark-Kernel-Univ-Central-Extension}
Using Proposition~\ref{Proposition-Universal-Central-Extension}, one sees that, for a perfect object $Y$, $H_{2}Y$ may be equivalently defined as the kernel of the universal central extension of~$Y$.
\end{remark}

Let $A$ be an abelian group object in $\Ac$. Recall that the sum $a+b$ of two elements $a,b\colon X\to A$ of a group $\hom (X,A)$ is the composite $m\comp (a,b)\colon {X\to A}$ of $(a,b)\colon X\to A\times A$ with the multiplication $m\colon A\times A\to A$ of $A$. Homming into $A$ defines a functor $\hom (\cdot,A)\colon \Ac^{\op}\to \Ab $. Given a simplicial object $S$ in $\Ac$, its image $\hom (S,A)$ is a cosimplicial object of abelian groups; as such, it has cohomology groups $H^{n}\hom (S,A)$.

\begin{definition}\label{Definition-Cotriple-Cohomology}
Let $\Ac$ be a semi-abelian category, monadic over $\Set$, and let $\Gb$ be the induced comonad. Let $X$ be an object of $\Ac$ and $A$ an abelian object. Consider $n\in \N_0$. We say that
\[
H^{n} (X,A)=H^{n-1}\hom (\ab (\Gb X),A)
\]
is the \emph{$n$-th cohomology group of $X$ with coefficients in $A$ (relative to the cotriple $\Gb$)}. This defines a functor $H^{n} (\cdot,A)\colon \Ac \to \Ab $, for any $n\in \N_0$. When it is clear which abelian group object $A$ is meant, we shall denote it just $H^{n} (\cdot)$.
\end{definition}

\begin{remark}\label{Remark-Is-Barr-Beck}
This is an instance of Barr and Beck's general definition of cotriple cohomology \cite{Barr-Beck}: $H^{n} (X,A)$ is nothing but the $n$-th cohomology group of $X$, with coefficients in the functor $\hom (\ab (\cdot),A)\colon {\Ac^{\op}\to \Ab}$, relative to the cotriple $\Gb$.
\end{remark}

\begin{proposition}\label{Proposition-First-Cohomology}
For any object $X$ of $\Ac$,
\[
H^{1} (X,A)\cong\hom (H_{1} X, A)\cong\hom (\ab (X),A)\cong\hom (X,A).
\]
If $X$ is projective then $H^{n} X=0$, for any $n\geq 2$.
\end{proposition}
\begin{proof}
The first isomorphism is a consequence of the fact that $\hom (\cdot, A)$ turns coequalizers in $\Ab(\Ac)$ into equalizers in $\Ab$. The second isomorphism follows from Theorem~\ref{Theorem-Homology} and the third one by adjointness of the functor $\ab$.

The second statement follows because if $X$ is projective then $\Gb X$ is contractible (see \cite{Barr-Beck}).
\end{proof}

The following result extends Theorem $12$ in {\cite{Carrasco-Homology}} and Theorem $1$ in \cite{AL}:
\begin{theorem}[Hochschild-Serre Sequence]\label{Theorem-Cohomology-Sequence}
Let
\[
\xymatrix{ 0 \ar[r] & K \ar@{{ >}->}[r]^-{k} & X \ar@{-{ >>}}[r]^-{f} & Y \ar[r] & 0 }
\]
be a short exact sequence in $\Ac$. There exists an exact sequence of abelian groups
\[
 0 \to<250> H^{1} Y \to/{{ >}->}/<250>^{H^{1}f} H^{1} X \to<250> \hom \bigl(\textstyle{\frac{K}{[K,X]}},A\bigr) \to<250> H^{2} Y \to<250>^{H^{2}f} H^{2} X
\]
 that depends naturally on the given short exact sequence.
\end{theorem}
\begin{proof}
The sequence (\ref{Sequence-simpB}) is degreewise split exact; hence homming into $A$ yields an exact sequence of abelian cosimplicial groups 
\[
 0 \to \hom (\ab\Gb Y,A)  \to/{{ >}->}/ \hom (\ab\Gb X, A) \to/{-{ >>}}/ \hom (K[\ab \Gb f],A) \to 0.
\]
This gives rise to an exact cohomology sequence
\[
0 \to H^{1} Y \to/{{ >}->}/^{H^{1}f} H^{1} X \to H^{0}\hom (K[\ab \Gb f],A) \to H^{2} Y \to^{H^{2}f} H^{2} X.
\]
By Proposition~\ref{Proposition-H_0-Sequence}, 
\[
H^{0}\hom (K[\ab \Gb f],A)\cong\hom (H_{0}K[\ab \Gb f], A)\cong\hom \bigl(\textstyle{\frac{K}{[K,X]}}, A\bigr),
\]
and the result follows.
\end{proof}

As a special case we get the following cohomological version of Hopf's formula.

\begin{corollary}\label{Corollary-CoHopf}
Let
\[
\xymatrix{ 0 \ar[r] & K \ar@{{ >}->}[r]^-{k} & X \ar@{-{ >>}}[r]^-{f} & Y \ar[r] & 0 }
\]
be a short exact sequence in $\Ac$, with $X$ a projective object. Then the sequence 
\[
\hom (X,A) \to<250> \hom \bigl(\textstyle{\frac{K}{[K,X]}},A\bigr) \to/{-{ >>}}/<250> H^{2} (Y,A) \to<250> 0
\]
is exact.
\end{corollary}
\begin{proof}
This follows immediately from the sequence in Theorem~\ref{Theorem-Cohomology-Sequence}, if we use Proposition~\ref{Proposition-First-Cohomology} which asserts that $H^{2} X=0$ when $X$ is projective.
\end{proof}

This means that an element of $H^{2} (Y,A)$ may be considered as an equivalence class $[a]$ of morphisms $a\colon {K}/{[K,X]}\to A$, where $[a]=[0]$ if and only if $a$ extends to $X$.

\section{The second cohomology group}\label{Section-Second-Cohomology}
In this section we characterize the second cohomology group $H^{2}(Y,A)$ of a group $Y$ with coefficients in an abelian object $A$ as the group $\CentrExt(Y,A)$ of isomorphism classes of central extensions of $Y$ by $A$. 

\begin{proposition}\label{Proposition-Baer-Sum}
Let $\Ac$ be a semi-abelian category and $Y$ an object of $\Ac$. Mapping an abelian object $A$ in $\Ac$ to the set $\CentrExt(Y,A)$ of isomorphism classes of central extensions of $Y$ by $A$ gives a finite product-preserving functor
\[
\CentrExt(Y,\cdot)\colon {\Ab(\Ac) \to \Set}.
\]
\end{proposition}
\begin{proof}
The functoriality of $\CentrExt(Y,\cdot)$ follows from Corollary~\ref{Corollary-Central-Extension-Pushout} and the Short Five Lemma: $\CentrExt(Y,1_{A})=1_{\CentrExt(Y,A)}$ is obvious, and $\CentrExt(Y,b\comp a)$ is equal to $\CentrExt(Y,b)\comp \CentrExt(Y,a)$ because the diagram with exact rows
\[
\xymatrix{0 \ar[r] & K\oplus B \ar[d]_-{a\oplus 1_{B}} \ar@{{ >}->}[r]^-{k\times 1_{B}} & X\times B \ar@{-{ >>}}[r] \ar[d]^-{\overline{a}\times 1_{B}} & Y \ar[r] \ar@{=}[d] & 0\\
0 \ar[r] & A\oplus B \ar@{{ >}->}[r]^-{k_{a}\times 1_{B}} \ar[d]_-{[b,1_{B}]} \ar@{}[rd]|>>{\pushout} & Z_{a}\times B \ar@{-{ >>}}[r] \ar[d]^-{\overline{\overline{ b}}} & Y \ar[r] \ar@{=}[d] & 0\\
0 \ar[r] & B \ar@{{ >}->}[r]_-{kk_{b}} & ZZ_{b} \ar@{-{ >>}}[r]_-{pp_{b}} & Y \ar[r] & 0}
\]
commutes, which yields a map from the induced pushout $Z_{b\comp a}$ to $ZZ_{b}$. This map is the needed isomorphism of central extensions. 

It is clear that $\CentrExt(Y,\cdot)$ preserves the terminal object. It also preserves binary products: the inverse of the map
\[
(\CentrExt(Y,\pi_{A}),\CentrExt(Y,\pi_{B}))\colon {\CentrExt(Y,A\times B)\to \CentrExt(Y,A)\times\CentrExt(Y,B)}
\]
is defined as follows. Given two central extensions
\[
0 \to A \to/{ >}->/^{k} X \to/-{ >>}/^{f} Y \to 0
\]
and
\[
0 \to B \to/{ >}->/^{l} Z \to/-{ >>}/^{g} Y \to 0,
\]
pulling back $f\times g$ along the diagonal $\Delta_{Y}=(1_{Y},1_{Y}) \colon  Y \to Y \times Y$ yields the diagram with exact rows
\[
\xymatrix{
0 \ar[r] & A\times B \ar@{{ >}->}[r]^-{k\times l} \ar@{=}[d] & X\times_{Y} Z \ar@{-{ >>}}[r]^-{h} \ar@{}[rd]|<<{\pullback} \ar[d]_{(\pi_{1},\pi_{2})} & Y \ar[r] \ar[d]^-{\Delta_{Y}} & 0\\
0 \ar[r] & A\times B \ar@{{ >}->}[r]_-{k\times l} & X\times Z \ar@{-{ >>}}[r]_-{f\times g} & Y\times Y \ar[r] & 0.}
\]
Let us denote the isomorphism class of a central extension $e$ as $\{e\}$. Then the couple $(\{f\},\{g\})$ is mapped to the isomorphism class $\{h\}$ of the map $h$, which is central as a pullback of the central extension $f\times g$. 
\end{proof}

The old definition of ``Baer sum'' \cite{Gerstenhaber} now becomes a simple instance of a general categorical fact: a finite product-preserving functor from an additive category to the category of sets factors uniquely over the category of abelian groups. This gives

\begin{proposition}
The functor $\CentrExt(Y,\cdot)$ factors uniquely over the forgetful functor ${\Ab\to \Set}$ to a functor ${\Ab(\Ac) \to \Ab}$, also denoted $\CentrExt(Y,\cdot)$.\noproof
\end{proposition}

Let us now explicitly describe the group structure on $\CentrExt(Y,K)$. Here $K$ is an abelian object; as such, it carries a multiplication $m\colon {K\times K\to K}$, which induces the map
\[
\CentrExt(Y,m)\comp (\CentrExt(Y,\pi_{1}),\CentrExt(Y,\pi_{2}))^{-1}
\]
as the multiplication (or rather, ``addition'') on $\CentrExt(Y,K)$. Let $f$ and $g$ be two central extensions as in the proof above, where now $A=B=K$, so that we may form the following diagram:
\[
\xymatrix{& 0 \ar[d] & 0 \ar[d]\\
& K[m] \ar@{=}[r] \ar@{{ >}->}[d]_{i=\ker m} & K[m] \ar@{{ >}->}[d]^{(k\times l)\comp i} \ar[r] & 0 \ar[d]\\
0 \ar[r] & K\times K \ar@{{ >}->}[r]^-{k\times l} \ar@{}[rd]|>>{\pushout} \ar@{-{ >>}}[d]_{m} & X\times_{Y}Z \ar@{.{ >>}}[d]^-{n} \ar@{-{ >>}}[r]^-{h} & Y \ar[r] \ar@{:}[d] & 0\\
0\ar@{.>}[r] & K \ar@{{ >}.>}[r]_{\overline{k\times l}} \ar[d] & W \ar@{.>}[d] \ar@{.{ >>}}[r]_-{f+g} & Y \ar@{.>}[r] & 0.\\
& 0 & 0}
\]
The arrow $(k\times l)\comp i$ is a kernel, thanks to Proposition~\ref{Proposition-Subobject-of-Central-is-Normal} and the fact that $h$ is central. The argument given in the proof of Proposition~\ref{Proposition-Reg-Epi-Central-Extension} shows that the bottom sequence is a central extension; its isomorphism class, denoted $\{f \} + \{ g \}$, clearly is the sum of the equivalence classes $\{f\}$ and $\{g\}$. As an immediate generalization of the case of groups, $\{f \} + \{ g \}$ could be called the \emph{Baer sum} of $\{f\}$ and $\{g\}$. (See Gerstenhaber \cite{Gerstenhaber} and, in a more general context, \cite{Bourn1999}.)

In summary: the sum of two classes $\{f \}$ and $\{g \}$ is the isomorphism class of the cokernel $f+g$ of the pushout $\overline{k\times l}$ of the arrow $k\times l$ along the multiplication $m$ of $K$.

\begin{theorem}\label{Theorem-H^2-Central-Extensions}
Let $\Ac$ be a semi-abelian category, monadic over $\Set$. Then the functor $H^{2} (Y,\cdot)$ is isomorphic to $\CentrExt (Y,\cdot)$.
\end{theorem}
\begin{proof}
We only have to prove that they are isomorphic as $\Set$-valued functors. To do so, let $A$ be an abelian object in $\Ac$, and $f\colon X\to Y$ a presentation of $Y$ with kernel $K$. Consider the reflection 
\[
0\to \tfrac{K}{[K,X]} \to/{ >}->/^{k} \tfrac{X}{[K,X]} \to/{-{ >>}}/^{\centr f} Y \to 0
\]
of $f$ into $\CentrExt (Y)$. In view of Corollary~\ref{Corollary-CoHopf}, we must show that there is a bijection $F$ from the set of equivalence classes $[a]$ of morphisms $a\colon {{K}/{[K,X]}\to A}$, where $[a]=[0]$ if and only if $a$ extends to $X$, to the set of isomorphism classes of central extensions of $Y$ by $A$.

The function $F$ is defined using Corollary~\ref{Corollary-Central-Extension-Pushout}: as $\centr f$ is a central extension, a morphism $a\colon {K}/{[K,X]}\to A$ gives rise to a central extension of $Y$ by $A$---of which the isomorphism class $F ([a])$ is the image of $[a]$ through $F$. 

$F$ is well-defined: if $[a]=[0]$ then $F ([a])=F ([0])$. Indeed, it is easily seen that $F ([0])$ is the isomorphism class of the central extension $\pi_{Y}\colon Y\times A\to Y$. If $a\colon K/[K,X]\to A$ factors over $X$ then it factors over $X/[K,X]$, and as a consequence the extension associated with $a$ has a split monic kernel. It follows that this extension is isomorphic to $\pi_{Y}\colon Y\times A\to Y$.

Finally, $F$ is a surjection because $X$ is projective and $\centr f$ is the reflection of $f$ into $\CentrExt (Y)$, and $F$ is injective because $F ([a])=\{\pi_{Y}\colon Y\times A\to Y\}$ entails that $a$ factors over $X$.
\end{proof}

\section{A universal coefficient theorem for cohomology in semi-abelian categories}\label{Section-Universal-Coefficient-Theorem}

This section treats the relationship between homology and cohomology.

Given two abelian objects $A$ and $C$ in $\Ac$, let $\ext (C,A)$ be the subgroup of $\CentrExt (C,A)$ $(\cong H^{2} (C,A)$) determined by the 
(isomorphism classes of) the extensions
\[
0 \to A  \to/{{ >}->}/ B \to/{-{ >>}}/ C \to 0
\]
of $C$ by $A$ lying in $\Ab (\Ac)$ (i.e.\ having the property that also $B$ belongs to $\Ab (\Ac)$).

Since the regular epimorphisms in $\Ab(\Ac) $ are just the regular epimorphisms of $\Ac$ that happen to lie in $\Ab(\Ac)$, the reflection $\ab (X)$ of a projective object $X$ of $\Ac$ is projective in $\Ab(\Ac)$. It follows that $\Ab(\Ac)$ has enough projectives if $\Ac$ has, and one may then choose a presentation
\begin{equation}\label{Abelian-Presentation}
0 \to R  \to/{{ >}->}/ F \to/{-{ >>}}/^{p} C \to 0
\end{equation}
of an abelian object $C$ in $\Ab(\Ac)$ instead of in $\Ac$.

\begin{proposition}\label{Proposition-Characterization-ext}
If $A$ is an abelian object and (\ref{Abelian-Presentation}) is a presentation in $\Ab (\Ac)$ of an abelian object $C$, then the sequence
\[
\hom (F,A) \to \hom (R,A) \to/{-{ >>}}/ \ext (C,A) \to 0
\]
is exact.
\end{proposition}
\begin{proof}
This is an application of Corollary~\ref{Corollary-Central-Extension-Pushout}. It suffices to note that the arrow $p$ is central, and that in a square
\[
\xymatrix{R \ar[d]_-{a} \ar@{{ >}->}[r] & F \ar[d] \\
A \ar@{{ >}->}[r] & Z_{a}}
\]
induced by Corollary~\ref{Corollary-Central-Extension-Pushout}, all objects are abelian. Indeed, $Z_{a}$ being an abelian object follows from the fact that $\Ab (\Ac)$ is closed under products and regular quotients in $\Ac$, and that $a$ can be decomposed as $[a,1_{A}]\comp l_{R}$.
\end{proof}

\begin{theorem}\label{Theorem-Universal-Coefficients}
If $Y$ is an object of $\Ac$ and $A$ is abelian then the sequence
\[
0 \to \ext (H_{1}Y,A)  \to/{{ >}->}/ H^{2} (Y, A) \to \hom (H_{2}Y,A)
\]
is exact.
\end{theorem}
\begin{proof}
One diagram says it all:
\[
\xymatrix{& \hom (H_{1} X,A) \ar@{=}[r] \ar[d] & \hom (X,A) \ar[d] \\
0 \ar[r] & \hom \bigl(\tfrac{K}{K\cap [X,X]},A\bigr) \ar@{{ >}->}[r] \ar@{-{ >>}}[d] & \hom \bigl(\tfrac{K}{[K,X]},A\bigr) \ar@{-{ >>}}[d] \ar[r] & \hom (H_{2}Y,A) \ar@{=}[d]\\
0\ar@{.>}[r] & \ext (H_{1}Y,A) \ar@{{ >}.>}[r] \ar[d] & H^{2} (Y,A) \ar@{.>}[r] \ar[d] & \hom (H_{2}Y,A).\\
& 0 & 0}
\]
Here $f\colon X\to Y$ is a presentation with kernel $K$, and the vertical sequences are exact by Proposition~\ref{Proposition-Characterization-ext}---the sequence
\[
0 \to \tfrac{K}{K\cap [X,X]}  \to/{{ >}->}/ H_{1}X \to/{-{ >>}}/ H_{1}Y \to 0
\]
being a presentation of $H_{1}Y$---and Corollary~\ref{Corollary-CoHopf}, respectively. The middle horizontal sequence is exact by the Hopf formula (Theorem~\ref{Theorem-Homology}) and the fact that, by the First Noether Isomorphism Theorem \cite[Theorem 4.3.10]{Borceux-Bourn},
\begin{equation}\label{Exact-Sequence-Sometimes-Split}
0 \to \tfrac{K\cap [X,X]}{[K,X]} \to/{{ >}->}/ \tfrac{K}{[K,X]} \to/{-{ >>}}/ \tfrac{K}{K\cap [X,X]} \to 0
\end{equation}
is an exact sequence in $\Ab(\Ac)$.  
\end{proof}

Recalling that an object $Y$ is perfect if and only if $H_{1}Y=0$, Theorem~\ref{Theorem-Universal-Coefficients} yields the following classical result.

\begin{corollary}\label{Corollary-Perfect-Homology-vs-Cohomology}
If $Y$ is a perfect object and $A$ is abelian then $H^{2} (Y,A)\cong\hom (H_{2}Y,A)$.
\end{corollary}
\begin{proof}
Comparing the Stallings-Stammbach sequence (\ref{Stallings-Stammbach-Sequence}) with Sequence (\ref{Exact-Sequence-Sometimes-Split}) and using that $Y$ is perfect we see that $K/ (K\cap [X,X])$ is isomorphic to $\ab (X)$. The latter object being projective in $\Ab(\Ac)$, the sequence (\ref{Exact-Sequence-Sometimes-Split}) is split exact in the abelian category $\Ab(\Ac)$. It follows that 
\[
0 \to \ext (H_{1}Y,A)  \to/{{ >}->}/ H^{2} (Y, A) \to/-{ >>}/ \hom (H_{2}Y,A)\to 0
\]
is a split exact sequence; but because $Y$ is perfect, $\ext (H_{1}Y,A)$ is zero.
\end{proof}

\begin{remark}\label{Remark-All-Central-Extensions}
In Section~\ref{Section-Second-Cohomology} we showed that $\CentrExt (Y,A)\cong H^{2} (Y,A)$. Combining this with Corollary~\ref{Corollary-Perfect-Homology-vs-Cohomology} one deduces that the category $\CentrExt (Y)$ of central extensions of a perfect object $Y$ is equivalent to the comma category $(H_{2}Y\downarrow \Ab (\Ac))$. This equivalence essentially follows from the universal property of the universal central extension of $Y$. As explained in Section~\ref{Section-Perfect-Case}, \emph{for any object} $Y$, the category $\CentrExt (Y)$ can be described as a category of discrete fibrations on the Galois groupoid of any presentation of $Y$, see \cite[Section $6$]{Janelidze-Kelly}.
\end{remark}

Given any object $Y$ and a presentation $F\to H_{1}Y$ with kernel $R$ in $\Ab(\Ac)$, Proposition~\ref{Proposition-Characterization-ext} entails the exactness of the sequence 
\[
\hom (F,A) \to \hom (R, A) \to/{-{ >>}}/ \ext (H_{1}Y,A) \to 0.
\]
If now, for every abelian object $A$ of $\Ac$, $\ext (H_{1}Y,A)$ is zero, then all functions ${\hom (F,A)\to \hom (R,A)}$ are surjections, which means that $R\to F$ is a split monomorphism. In this case it follows that $F=R\oplus H_{1}Y$ and $H_{1}Y$ is projective in $\Ab(\Ac)$. As a consequence, we get the following partial converse to Corollary~\ref{Corollary-Perfect-Homology-vs-Cohomology}.

\begin{corollary}\label{Corollary-Perfect-Homology-vs-Cohomology-2}
If, for every abelian object $A$ of $\Ac$, $H^{2} (Y,A)\cong\hom (H_{2}Y,A)$, then $H_{1}Y$ is projective in $\Ab(\Ac)$.\noproof 
\end{corollary}

\small

\providecommand{\bysame}{\leavevmode\hbox to3em{\hrulefill}\thinspace}
\providecommand{\MR}{\relax\ifhmode\unskip\space\fi MR }
\providecommand{\MRhref}[2]{%
  \href{http://www.ams.org/mathscinet-getitem?mr=#1}{#2}
}
\providecommand{\href}[2]{#2}

\end{document}